\documentstyle{article} 
\begin{document}
\title{On the section of a cone} 
\author{Alejandro Rivero}
\maketitle
\begin{abstract}
A problem from Democritus is used to illustrate the
building and use of infinitesimal covectors. 
\end{abstract}

The Friday before Passover I was forced to make some bureaucratic consultations 
in our Ministery of Defence. So I landed Tuesday in our desolated airport at 
Zaragoza and, two days after, I took the train to Madrid, hoping to get 
lodgement in the house of Miss Ana Leal in the folkish town of {\it La Latina}.

Such happenings use to be intellectually exciting, and this one was not 
exception. Miss Leal suggested extending some hours  the visit in order to be 
able to attend a lecture of Agustin  Garcia Calvo in Lavapies. This Agustin is a 
well-known classical linguist, and a kind of anarchist philosopher, who likes to 
teach in the Greek stile, and we had already enjoyed his {\it tertulias} in the 
old institution of the {\it Ateneo de Madrid}. This one was supposed to be a 
more technical lecture, addressed to secondary school teachers.

Indeed, it was a very amenable lecture. According to my notebook, he made a good 
point of the use of language for political control, opposing the vocabulary 
against grammar, the former owned by the power to construct the reality, the 
latter unconsciently managed by the people driving a "raison en marche". Being 
practising mathematician, one can easily to feel this confrontation; gramatician 
placeholders, names, adjectives, etc, are not very different of our variables 
and constants, and our whole fight  is to leave all the weight over proofs, over 
our grammar, avoiding to get any conclusions from the loose vocabulary of 
definitions. I meditated on these parallelisms while  hearing the linguist' 
admonitions. 

Then, Agustin centred in the lecture main theme, the teaching of philosophy in 
secondary education, and somehow malevously suggested three examples to be 
proposed to students. The Lewis Carroll approach to Zeno paradoxes, the Zeno 
paradox itself (beautifully glossed as {\it no se vive mientras se besa, no se 
besa mientras se vive}, one does not live while kissing, one does not kiss while 
living) and, for my surprise, the dilemma of the cone from Democritus!

Heath, following Plutarch, enunciates it asking what happens if we cut
a cone using a plane parallel and very close to the basis. Is the 
resulting circle equal to the one of the basis? Our XXth century
presocratic philosopher prefers a more intrincate set up; let me to
 go to my notes and remember this. Take a cone and cut it through
a plane, which for 
simplicity we can still take parallel to the basis. Now, look at
 the resulting 
figures, a smaller cone with basis B and a conical trunk with top B'. The 
question is, are the circles B and B' equal or unequal?

In any case, if they are not equal, it results that some discontinuity happens in the 
complete, joined, cone, and the generating line  should present jumps, small 
steps. But if, on the contrary, both surfaces are not unequal, their fusion 
should build a cylinder, no a cone.

To solve this paradox, we can negate the possibility of the described action. We 
can claim that the cone is a real figure,  and then it is not proofed that it is 
a mathematical cone, and its generating line could really be irregular. But  
the mathematical problem still exists, and we can ask about the ideal 
cone\footnote{We can claim that 
 "we can see it with the eye of the mind; and we know, by 
force of demonstration, that it cannot be otherwise", as Democritus
himself claimed for the tangent of the circle.}. 

Again, in this setup, we can negate the starting point and claim that the cut  is really the 
intersection of the plane with the cone. There are no two  surfaces to be 
compared.

But then the paradox can be got again by using a Gedankenexperiment.  Instead of 
Plutarch quote, let me get the same music, if not the notes, from the last point 
of Agustin: Imagine we cut a carrot, or a turnip with a cutter, so we can not 
deny we have two surfaces. In principle the cutter retires some slice of
matter from the carrot, slice thickness being related to the one of the
cutter. Imagine we make the cutter thinner and thinner, so no 
mass is moved out of the carrot when we cut it. Now imagine the same operation 
on the mathematical setup, we have two surfaces and the paradox again.

And from here we are in our own\footnote{I intend to hide some
complex or distracting comments under the carpet of the footnotes. The reader
could prefer to avoid them in a first reading}. 

The result of the progressive thinning of the cutter has been a pair of planes 
becoming progressively nearer. This is to be noted: the operative problem 
involves no a plane, but two planes approaching one to other. 

This figure, a pair of parallel planes, is known to mathematics following 
Shouten and Golab as being a covector (in a tridimensional space). Technically, 
it is specified by giving a unit segment (axial vector) perpendicular to the 
plane, and then a modulus measuring the separation of both planes and a support 
point where the first of the planes lies\footnote{The specification of the 
axial vector changes peculiarly when we make a change of coordinates of the 
system, fitting the usual definition of covectors. And the space of covectors 
(n-1 segments figures as specified) is dual to the space of vectors (1 
dimensional oriented segments). An equivalence relationship can be added to get 
a space of free covectors, but here this step is not needed. We can say 
that two covectors can be added when the final plane of one coincides 
with the starting plane of the other, then fusing it to make a grosser cutter.
To be honest, this restriction is stronger than the usual for 
"free" covectors, and in fact it reflect that we are interested in a slightly 
looser structure, which we could call q-covectors (the q making reference to a 
scale of the thickness of the cutter, and -indirectly- to the deformed 
differential calculus of Majid)}.

Now, lets to make the cutter slimmer and slimmer. Then the modulus of our 
covectors goes to zero, but we still have two planes.

This is in fact the resolution of the paradox. We distinguish  between a 
cylinder and a cone because we have more information: the continuity should be 
claimed between one surface of the pair and the next one of the following 
"slice". The difference between a cone and a cylinder resides in the internal 
structure of the pair. When the pair becomes so-to-say infinitesimal, both 
planes of the slice live over the same points of the space. All the 
infinitesimal slices over the same area can be added without taking care of the 
fusion condition above, and then get additional structure\footnote{ 
Mathematicians call this space of infinitesimal covectors  over a 
point "cotangent space", the whole set being the "cotangent bundle". An
application selecting one covector over each point is called a "differential 
form".}.

Of course, it must be seen that there is something in the limit of approaching 
surfaces, i.e., we must to give sense to this limit and proof there are really 
something over a single cut. In our modern XXth century we could jump
directly to use scaling 
transformations in the spirit of Wilson-Kogut. But perhaps it is 
better to start from Archimedian methods, which the reader can  enjoy in the 
interesting book of T.L. Heath. For instance, we can see how the volume of a 
conoid can be extracted by using our bifacial knife. And, as we ignore the full  
detail of Greek methods\footnote{Heath quotes Wallis  regretting that "nearly 
all the ancients so hid from posterity their method of Analysis (though it is 
clear that they had one) that more modern mathematicians found it easier to 
invent a new  Analysis than to seek out the old". Indeed, the lack of texts is 
surprising, all a branch of reason cleared white as the recycled folia where 
"The Method" was found in 1909, {\it palimpsesta sunt, scriptura antiqua 
(litteris minusculis s. X) aqua tantum diluta plerumque oculis intentis discipi  
potest (de foll 1..., 119-122 tamen desperandum mihi erat)}, a inmense cleaning 
which justifies Wallis' paranoia. But again, even Newton kept secret his own 
method, until that Leibnitz developments forced him to show it.}, it could be 
perhaps forgiven if we avoid refilling the discussion with Leibnitzian meat, 
trying instead to keep the spicy flavour of our local cooking.

Imagine again the cone, divided into slices of finite size, lets say using our 
finite cutter. Each slice can be fitted between two cylindrical ones, a smaller 
one, with takes as basis the small circle of the conic slice, and a greater one, 
taking as basis the big circle.

By joining the cylinders we have two circular ''ziggurats'', a small one inscribed 
inside the cone, and a greater one circumscribed on it. We can consider then the 
calculation of a quantity such all the volume of these whole figures from the 
corresponding quantity of the pieces.

The difference between the circumscribed and the inscribed figure amounts only 
to the greater slice of the circumscribed one. This is because each cylinder of 
the inscribed one is equal to the previous one of the circumscribed one.
Here we see that the importance of the correct pasting condition: it must be 
between the circle of one slice and the immediate of the following one. Only in 
this manner the subtraction keeps control, all the difference being the volume of
only one slice. Thus when the cutter thickness goes to
zero, so goes the difference between figures, and their volume converges to
the volume of the cone.

Lets examine this convergence with more detail. It involves two operations: 
to increase the number of slices, to decrease the width of the slices. Both
operations related, of course, because the product is the height of the
cone. But here we do not see the structure of the limiting objects, so
it is possible yet to hold some doubts about the process. An alternative
approach is the averaging method of Wilson \footnote{The interested
reader can see some examples in the article published by Wilson 
himself in 1979 in the Scientific American}. Two consecutive
cylinders can be substituted by an unique cylinder averaging them,
i.e., with a volume that is the sum of the volumes of both cilinders and
a thickness equal to the union of them, thus double of the original
one. 

Aplying this procedure to the whole "ziggurat" we get a new figure
which is no more inscribed (or circunscribed) to the cone, but
has the same volume that the starting one. 

Now, this method can be used to control the limit process in the
following way: we choose an arbitrary scale of thickness, say
for instance the one half of the height of the cone, and for 
each "ziggurat" in the converging series we apply the averaging
until we get back to a figure composed of cylinders of the
choosen thickness.

This new series of figures\footnote{The new series could be
called "renormalized", if we call the original "bare"} is composed
of finite objects, each one having equal number of cylinders
and cylinder thickness being the same in every figure. Then
the limit process of this series is not affected by the two
infinities, in thickness and in number of cylinders, that were
incrusted into the previous series. 
 Even if we do not believe in the infinitesimal slices, we should
have not problem admiting the regularized slices built at given, but
arbitrary, scale. 

Readers could note that a "ziggurat" composed of cylinders of equal
radius should be invariant under the Wilson transformation, only
the external, nominal, scale of thickness changes. A deeper 
examination would show that the existence of this
set of invariant figures is the key for the convergence of
the whole process\footnote{By iteration of the transformation,
we finish over some invariant figure, and the slice we are 
slimming is also similar to the invariant figure. If the transformation
were given by a continous group, we should see trajectories on
the space the figures, approaching the trajectory of invariant
cylinders, and the renormalized series would be a line cutting
across trajectories and converging to a point in the invariant
trajectory. In our case, with discrete transformations, we can still
hit into the invariant trajectory by choosing the length (for
instance, using as fundamental length the heigth of the cone,
instead of one half as above).
}. In some sense, this invariant shape is the
amplification, to a finite scale, of the infinitesimal cylinders
of the first convergence process. We can choose
the arbitrary reference scale as near to zero as we wish, so
its limit zero\footnote{Which we could be tempted to
call "classical limit" instead of the usual "continous limit"} 
can be interpreted as the home of such "differentials". In fact
the existence of limit for our finite series depends of the
existence of the line of invariant figures, and the existence
of this line relates to the existence of a fixed point, from 
which the line starts. Such fixed point can be linked the above
suggested zero limit. All the pieces of the puzzle fit together.

Note that equal that we can cut the cone, we can also cut its axis, the only 
difference being that the former lives in three dimensions, the latter in one, 
then the covectors over this latter are specified by pairs of points instead of 
planes. Following conventions, we can call dV to the infinitesimal slice over 
the cone,  and dz to the slice over the axis. It is also usual to write $ dV= 
{\partial V \over \partial z} dz = A(z) dz$ but of course such writing must also be 
justified.  

Another stroke could be draw if we choose to imagine the cone as developing in 
the time, i.e., the axis being some kind of temporal direction. Thus each cut is 
a circle which grows in the time, and the connection with Zeno paradox becomes 
evident \footnote{Sketching a parallelism with modern physics terminology, we 
could say that Democritus paradox is the "Wick-rotated" version of this from 
Zeno. The axis should be the "imaginary time", and volume and area should 
perhaps correspond to position and velocity.}

And further developments could be done, for instance connecting the scaling 
procedure to the ones currently in use in theoretical physics, or to work out 
the q-covectors composition rule aiming to a cotangent groupoid similar to the 
tangent groupoid raised by Monsieur Alain Connes. The chain of reasonments
is enough tight to rule out impossible relations and, as an old friend used
to say, when impossible is ruled out, the only remaining thing is the
answer. In this footing, we could follow up trying to build the
arguments in four
dimensional spaces and within field theory, from where we have already 
taken part of our terminology.

For sure that readers can imagine a lot of additional quests.

So, which our conclusion is? Well, we have seen how the discussion about
such an old problem becomes an argument for 
the teaching of modern mathematics and  physics. Perhaps this is the
whole point of this 
note, although surely it was not the one of Garcia Calvo when spelling this old 
tale to the philosophical audience. But again, mathematics works in its own 
pace, independently of our own intentions, just as  sometimes science gets to be 
taught independently of the intention of educational programs. Well, this one 
was probably the very point of Agustin, ours is only to give voice to the math 
through ourselves. He says {\it dejarse hablar}.

\end{document}